\documentclass[11pt]{article}%
\usepackage{amsfonts}
\usepackage{amssymb}
\usepackage{amsmath}
\usepackage{amscd}
\usepackage{latexsym}
\usepackage[doublespacing]{setspace}
\usepackage{theorem}
\usepackage{natbib}
\usepackage[colorlinks=true,urlcolor=black,linkcolor=black,pdfpagemode="None",pdfstartview=FitH]{hyperref}
\usepackage{graphicx}
\usepackage{caption}
\usepackage[font={normalsize},labelfont={normalsize},textfont={normalsize}]{subfig}%
\providecommand{\U}[1]{\protect\rule{.1in}{.1in}}
\newtheorem{theorem}{Theorem}

\newtheorem{lemma}{Lemma}

\newcommand{\I}{{\rm 1\hspace*{-0.4ex}\rule{0.1ex}{1.52ex}\hspace*{0.2ex}}}

\hypersetup{pdfborder = {0 0 0},colorlinks=true,linkcolor=blue,citecolor=blue}
\renewcommand{\cite}{\citet*}
\makeatletter
\renewcommand*{\numberline}[1]{\hb@xt@3.5em{#1\hfil\qquad}}
\makeatother
\begin{document}

\title{Bootstrap-Based Inference for Cube Root Asymptotics\thanks{A previous version
of this paper circulated under the title \textquotedblleft Bootstrap-Based
Inference for Cube Root Consistent Estimators\textquotedblright. For comments
and suggestions, we are grateful to Mehmet Caner, Andreas Hagemann, Kei
Hirano, Bo Honor\'{e}, Guido Imbens, Guido Kuersteiner, Mykhaylo Shkolnikov,
Ronnie Sircar, Ulrich M\"{u}ller, Whitney Newey, and participants at various
conferences, workshops, and seminars. Cattaneo gratefully acknowledges
financial support from the National Science Foundation through grants
SES-1459931 and SES-1947805, and Jansson gratefully acknowledges financial
support from the National Science Foundation through grants SES-1459967 and
SES-1947662 and the research support of CREATES (funded by the Danish National
Research Foundation under grant no. DNRF78).}}
\author{Matias D. Cattaneo\thanks{Department of Operations Research and Financial
Engineering, Princeton University.}
\and Michael Jansson\thanks{Department of Economics, University of California at
Berkeley and CREATES.}
\and Kenichi Nagasawa\thanks{Department of Economics, University of Warwick.}}
\maketitle

\begin{abstract}
This paper proposes a valid bootstrap-based distributional approximation for
$M$-estimators exhibiting a \cite{Chernoff_1964_AISM}-type limiting
distribution. For estimators of this kind, the standard nonparametric
bootstrap is inconsistent. The method proposed herein is based on the
nonparametric bootstrap, but restores consistency by altering the shape of the
criterion function defining the estimator whose distribution we seek to
approximate. This modification leads to a generic and easy-to-implement
resampling method for inference that is conceptually distinct from other
available distributional approximations. We illustrate the applicability of
our results with four examples in econometrics and machine learning.

\end{abstract}


%TCIMACRO{\TeXButton{setup}{\thispagestyle{empty}
%\setcounter{page}{0}
%\newpage\setlength{\abovedisplayskip}{5pt}
%\setlength{\belowdisplayskip}{5pt}}}%
%BeginExpansion
\thispagestyle{empty}
\setcounter{page}{0}
\newpage\setlength{\abovedisplayskip}{5pt}
\setlength{\belowdisplayskip}{5pt}%
%EndExpansion

\section{Introduction\label{[Section] Introduction}}

In a seminal paper, \cite{Kim-Pollard_1990_AoS} studied estimators exhibiting
\textquotedblleft cube root asymptotics\textquotedblright. These estimators
not only have a non-standard rate of convergence, but also have the property
that rather than being Gaussian their limiting distributions are of
\cite{Chernoff_1964_AISM} type; i.e., the non-Gaussian limiting distribution
is that of the maximizer of a Gaussian process. Kim and Pollard's results
cover not only celebrated examples such as maximum score estimator of
\cite{Manski_1975_JoE} and the isotonic density estimator of
\cite{Grenander_1956_SAJ}, but also more contemporary estimators arising in
examples related to classification problems in machine learning
%TCIMACRO{\TeXButton{\citep*{Mohammadi-vandeGeer_2005_JMLR}}{\citep
%*{Mohammadi-vandeGeer_2005_JMLR}}}%
%BeginExpansion
\citep*{Mohammadi-vandeGeer_2005_JMLR}%
%EndExpansion
, nonparametric inference under shape restrictions
%TCIMACRO{\TeXButton{\citep*{Groeneboom-Jongbloed_2018_SS}}{\citep
%*{Groeneboom-Jongbloed_2018_SS}}}%
%BeginExpansion
\citep*{Groeneboom-Jongbloed_2018_SS}%
%EndExpansion
, massive data $M$-estimation framework
%TCIMACRO{\TeXButton{\citep*{Shi-Lu-Song_2018_JASA}}{\citep
%*{Shi-Lu-Song_2018_JASA}}}%
%BeginExpansion
\citep*{Shi-Lu-Song_2018_JASA}%
%EndExpansion
, and maximum score estimation in high-dimensional settings
%TCIMACRO{\TeXButton{\citep*{Mukherjee-Banerjee-Ritov_2019_wp}}{\citep
%*{Mukherjee-Banerjee-Ritov_2019_wp}}}%
%BeginExpansion
\citep*{Mukherjee-Banerjee-Ritov_2019_wp}%
%EndExpansion
. Moreover, \cite{Seo-Otsu_2018_AoS} recently generalized
\cite{Kim-Pollard_1990_AoS} to allow for $n$-varying objective functions ($n$
denotes the sample size), further widening the applicability of cube-root-type
asymptotics. For example, their results cover the conditional maximum score
estimator of \cite{Honore-Kyriazidou_2000_ECMA}.

An important feature of Chernoff-type asymptotic distributional approximations
is that the covariance kernel of the Gaussian process characterizing the
limiting distribution often depends on an infinite-dimensional nuisance
parameter. From the perspective of inference, this feature of the limiting
distribution represents a nontrivial complication relative to the conventional
asymptotically normal case, where the limiting distribution is known up to the
value of a finite-dimensional nuisance parameter (namely, the covariance
matrix of the limiting distribution). The dependence of the limiting
distribution on an infinite-dimensional nuisance parameter implies that
resampling-based distributional approximations seem to offer the most
attractive approach to inference in estimation problems exhibiting cube root
asymptotics. Unfortunately, however, the standard nonparametric bootstrap is
well known to be invalid in this setting
%TCIMACRO{\TeXButton{\citep*[e.g.,][]{Bootstrap failure cites}}{\citep
%*{Abrevaya-Huang_2005_ECMA,Leger-MacGibbon_2006_CJS,Kosorok_2008_BookCh,Sen-Banerjee-Woodroofe_2010_AoS}%
%}}%
%BeginExpansion
\citep
*{Abrevaya-Huang_2005_ECMA,Leger-MacGibbon_2006_CJS,Kosorok_2008_BookCh,Sen-Banerjee-Woodroofe_2010_AoS}%
%EndExpansion
. The purpose of this paper is to propose a generic and easy-to-implement
bootstrap-based distributional approximation applicable in the context of cube
root asymptotics.

As does the familiar nonparametric bootstrap, the method proposed herein
employs bootstrap samples of size $n$ from the empirical distribution
function. But unlike the nonparametric bootstrap, which is inconsistent, our
method offers a consistent distributional approximation for estimators
exhibiting cube root asymptotics. Consistency is achieved by altering the
shape of the criterion function defining the estimator whose distribution we
seek to approximate. Heuristically, the method is designed to ensure that the
bootstrap version of a certain empirical process has a mean resembling the
large sample version of its population counterpart. The latter is quadratic in
the problems we study, and known up to the value of a certain matrix. As a
consequence, the only ingredient needed to implement the proposed
\textquotedblleft reshapement\textquotedblright\ of the objective function is
a consistent estimator of the unknown matrix entering the quadratic mean of
the empirical process. Such estimators turn out to be generically available
and easy to compute.

This paper is not the first to propose a consistent resampling-based
distributional approximation for cube-root-type estimators. For canonical cube
root asymptotic problems, the best known consistent alternative to the
nonparametric bootstrap is probably subsampling
%TCIMACRO{\TeXButton{\citep{Politis-Romano_1994_AoS}}{\citep
%{Politis-Romano_1994_AoS}}}%
%BeginExpansion
\citep{Politis-Romano_1994_AoS}%
%EndExpansion
, whose applicability was pointed out by
\cite{Delgado-RodriguezPoo-Wolf_2001_EL}. Related applicable methods are the
$m$ out of $n$ bootstrap
%TCIMACRO{\TeXButton{\citep{Bickel-Gotze-vanZwet_1997_SSinica}}{\citep
%{Bickel-Gotze-vanZwet_1997_SSinica}}}%
%BeginExpansion
\citep{Bickel-Gotze-vanZwet_1997_SSinica}%
%EndExpansion
, whose applicability was discussed and extended by \cite{Lee-Pun_2006_JASA}
and \cite{Lee-Yang_2020_AoS}, the rescaled bootstrap
%TCIMACRO{\TeXButton{\citep{Duembgen_1993_PTRF}}{\citep{Duembgen_1993_PTRF}}}%
%BeginExpansion
\citep{Duembgen_1993_PTRF}%
%EndExpansion
, and the numerical bootstrap
%TCIMACRO{\TeXButton{\citep{Hong-Li_2020_AoS}}{\citep{Hong-Li_2020_AoS}}}%
%BeginExpansion
\citep{Hong-Li_2020_AoS}%
%EndExpansion
. In addition, case-specific (smooth or non-standard) bootstrap methods have
been proposed for leading examples such as monotone density estimation
%TCIMACRO{\TeXButton{\citep
%{Kosorok_2008_BookCh,Sen-Banerjee-Woodroofe_2010_AoS}}{\citep
%{Kosorok_2008_BookCh,Sen-Banerjee-Woodroofe_2010_AoS}}}%
%BeginExpansion
\citep{Kosorok_2008_BookCh,Sen-Banerjee-Woodroofe_2010_AoS}%
%EndExpansion
, maximum score estimation
%TCIMACRO{\TeXButton{\citep{Patra-Seijo-Sen_2018_JoE}}{\citep
%{Patra-Seijo-Sen_2018_JoE}}}%
%BeginExpansion
\citep{Patra-Seijo-Sen_2018_JoE}%
%EndExpansion
, and the current status model
%TCIMACRO{\TeXButton{\citep{Groeneboom-Hendrickx_2018_AoS}}{\citep
%{Groeneboom-Hendrickx_2018_AoS}}}%
%BeginExpansion
\citep{Groeneboom-Hendrickx_2018_AoS}%
%EndExpansion
. For the more generic cube-root-type estimators analyzed in
\cite{Seo-Otsu_2018_AoS}, subsampling appears to be the only method available,
and indeed the authors discuss in their concluding remarks the need for (and
importance of) developing resampling methods based on the standard
nonparametric bootstrap. Our paper appears to be the first to provide one such method.

Like ours, each of the resampling methods mentioned above can be viewed as
offering a \textquotedblleft robust\textquotedblright\ alternative to the
standard nonparametric bootstrap but, unlike ours, existing methods achieve
consistency by modifying the distribution used to generate the bootstrap
sample. In contrast, our bootstrap-based method achieves consistency by means
of an analytic modification of the objective function used to construct the
bootstrap-based distributional approximation. As further discussed below, this
approach results in a procedure that is conceptually related to the bootstrap
methods developed by \cite{Andrews-Soares_2010_ECMA} and
\cite{Fang-Santos_2019_REStud} in other econometrics contexts.

Implementation of our procedure is not computationally demanding. Indeed, the
only ingredient needed to implement our modification on the objective function
is a consistent estimator of a certain Hessian matrix. We propose a generic
estimator based on numerical derivatives and present a consistency result as
well as an approximate mean square error expansion for that estimator. In
addition, we illustrate how example-specific features can be sometimes
exploited to construct alternative estimators.

The paper proceeds as follows. Section \ref{[Section] Heuristics} is heuristic
and outlines the main idea underlying our approach in the $M$-estimation
setting of \cite{Kim-Pollard_1990_AoS}. Section \ref{[Section] Main Result}
then makes that heuristic discussion rigorous in a more general setting
similar to that of \cite{Seo-Otsu_2018_AoS}. Section \ref{[Section] Examples}
illustrates our bootstrap-based inference method with four examples: the
maximum score estimator of
%TCIMACRO{\TeXButton{\cite{Manski_1975_JoE,Manski_1985_JoE}}{\cite
%{Manski_1975_JoE,Manski_1985_JoE}}}%
%BeginExpansion
\cite{Manski_1975_JoE,Manski_1985_JoE}%
%EndExpansion
, the conditional maximum score panel data estimator of
\cite{Manski_1987_ECMA}, the conditional maximum score dynamic panel data
estimator of \cite{Honore-Kyriazidou_2000_ECMA}, and the classification
estimator of \cite{Mohammadi-vandeGeer_2005_JMLR}. Section
\ref{[Section] Simulations} reports simulation evidence for the case of the
maximum score estimator, and Section \ref{[Section] Conclusion} concludes.
Section \ref{[Section] Proof of Theorem 1} describes the proof of our main
result, while the supplemental appendix contains omitted proofs and details.

\section{Heuristics\label{[Section] Heuristics}}

Suppose $\mathbf{\theta}_{0}\in\mathbf{\Theta}\subseteq\mathbb{R}^{d}$ is an
estimand admitting the characterization%
\begin{equation}
\mathbf{\theta}_{0}=\operatorname*{argmax}\limits_{\mathbf{\theta}%
\in\mathbf{\Theta}}M_{0}(\mathbf{\theta}),\qquad M_{0}(\mathbf{\theta
})=\mathbb{E}[m_{0}(\mathbf{z},\mathbf{\theta})],\label{Estimand}%
\end{equation}
where $m_{0}$ is a known function, and where $\mathbf{z}$ is a random vector
of which a random sample $\mathbf{z}_{1},\ldots,\mathbf{z}_{n}$ is available.
Studying estimation problems of this kind for non-smooth $m_{0},$
\cite{Kim-Pollard_1990_AoS} gave conditions under which the $M$-estimator%
\[
\mathbf{\hat{\theta}}_{n}=\operatorname*{argmax}\limits_{\mathbf{\theta
\in\Theta}}\hat{M}_{n}(\mathbf{\theta}),\qquad\hat{M}_{n}(\mathbf{\theta
})=\frac{1}{n}\sum_{i=1}^{n}m_{0}(\mathbf{z}_{i},\mathbf{\theta}),
\]
exhibits cube root asymptotics:%
\begin{equation}
\sqrt[3]{n}(\mathbf{\hat{\theta}}_{n}-\mathbf{\theta}_{0})\rightsquigarrow
\operatorname*{argmax}\limits_{\mathbf{s}\in\mathbb{R}^{d}}\{\mathcal{G}%
_{0}(\mathbf{s})+\mathcal{Q}_{0}(\mathbf{s})\},\label{Cube root asymptotics}%
\end{equation}
where $\rightsquigarrow$ denotes weak convergence, $\mathcal{G}_{0}$ is a
non-degenerate zero-mean Gaussian process, and $\mathcal{Q}_{0}(\mathbf{s}%
)=-\mathbf{s}^{\prime}\mathbf{H}_{0}\mathbf{s}/2,$ where $\mathbf{H}%
_{0}=-\partial^{2}M_{0}(\mathbf{\theta}_{0})/\partial\mathbf{\theta}%
\partial\mathbf{\theta}^{\prime}.$

Whereas the matrix $\mathbf{H}_{0}$ governing the shape of $\mathcal{Q}_{0}$
is finite-dimensional, the covariance kernel of $\mathcal{G}_{0}$ in
(\ref{Cube root asymptotics}) typically involves infinite-dimensional unknown
quantities. As a consequence, the limiting distribution of $\mathbf{\hat
{\theta}}_{n}$ tends to be more difficult to approximate than Gaussian
distributions, implying in turn that basing inference on $\mathbf{\hat{\theta
}}_{n}$ is more challenging under cube root asymptotics than in the more
familiar case where $\mathbf{\hat{\theta}}_{n}$ is $\sqrt{n}$-consistent and
asymptotically normally distributed.

As a candidate method of approximating the distribution of $\mathbf{\hat
{\theta}}_{n},$ consider the nonparametric bootstrap. To describe it, let
$\mathbf{z}_{1,n}^{\ast},\ldots,\mathbf{z}_{n,n}^{\ast}$ denote a random
sample from the empirical distribution of $\mathbf{z}_{1},\ldots\mathbf{z}%
_{n}$ and let the natural bootstrap analogue of $\mathbf{\hat{\theta}}_{n}$ be
denoted by%
\[
\mathbf{\hat{\theta}}_{n}^{\ast}=\operatorname*{argmax}\limits_{\mathbf{\theta
\in\Theta}}\hat{M}_{n}^{\ast}(\mathbf{\theta}),\qquad\hat{M}_{n}^{\ast
}(\mathbf{\theta})=\frac{1}{n}\sum_{i=1}^{n}m_{0}(\mathbf{z}_{i,n}^{\ast
},\mathbf{\theta}).
\]
Then, the nonparametric bootstrap estimator of $\mathbb{P}[\mathbf{\hat
{\theta}}_{n}-\mathbf{\theta}_{0}\leq\cdot]$ is given by $\mathbb{P}_{n}%
^{\ast}[\mathbf{\hat{\theta}}_{n}^{\ast}-\mathbf{\hat{\theta}}_{n}\leq\cdot],$
where $\mathbb{P}_{n}^{\ast}$ denotes a probability computed under the
bootstrap distribution conditional on the data. As is well documented,
however, this estimator is inconsistent under cube root asymptotics
%TCIMACRO{\TeXButton{\citep[e.g.,][]{Bootstrap failure cites}}{\citep
%*{Abrevaya-Huang_2005_ECMA,Leger-MacGibbon_2006_CJS,Kosorok_2008_BookCh,Sen-Banerjee-Woodroofe_2010_AoS}%
%}}%
%BeginExpansion
\citep
*{Abrevaya-Huang_2005_ECMA,Leger-MacGibbon_2006_CJS,Kosorok_2008_BookCh,Sen-Banerjee-Woodroofe_2010_AoS}%
%EndExpansion
.

For the purpose of giving a heuristic, yet constructive, explanation of the
inconsistency of the nonparametric bootstrap, it is helpful to recall that a
proof of (\ref{Cube root asymptotics}) can be based on the representation%
\begin{equation}
\sqrt[3]{n}(\mathbf{\hat{\theta}}_{n}-\mathbf{\theta}_{0}%
)=\operatorname*{argmax}\limits_{\mathbf{s}\in\mathbb{R}^{d}}\{\hat{G}%
_{n}(\mathbf{s})+Q_{n}(\mathbf{s})\},\label{Representation: Thetahat}%
\end{equation}
where, for $\mathbf{s}$ such that $\mathbf{\theta}_{0}+\mathbf{s}n^{-1/3}%
\in\mathbf{\Theta,}$%
\begin{equation}
\hat{G}_{n}(\mathbf{s})=n^{2/3}[\hat{M}_{n}(\mathbf{\theta}_{0}+\mathbf{s}%
n^{-1/3})-\hat{M}_{n}(\mathbf{\theta}_{0})-M_{0}(\mathbf{\theta}%
_{0}+\mathbf{s}n^{-1/3})+M_{0}(\mathbf{\theta}_{0})]\label{Ghat}%
\end{equation}
is a zero-mean random process, while%
\begin{equation}
Q_{n}(\mathbf{s})=n^{2/3}[M_{0}(\mathbf{\theta}_{0}+\mathbf{s}n^{-1/3}%
)-M_{0}(\mathbf{\theta}_{0})]\label{Q}%
\end{equation}
is a non-random function that is correctly centered in the sense that
$\operatorname*{argmax}_{\mathbf{s}\in\mathbb{R}^{d}}Q_{n}(\mathbf{s}%
)=\mathbf{0.}$ In cases where $m_{0}$ is non-smooth but $M_{0}$ is smooth,
$\hat{G}_{n}$ and $Q_{n}$ are usually asymptotically Gaussian and
asymptotically quadratic, respectively, in the sense that%
\begin{equation}
\hat{G}_{n}(\mathbf{s})\rightsquigarrow\mathcal{G}_{0}(\mathbf{s}%
)\label{Convergence: Ghat}%
\end{equation}
and%
\begin{equation}
Q_{n}(\mathbf{s})\rightarrow\mathcal{Q}_{0}(\mathbf{s}).\label{Convergence: Q}%
\end{equation}
Under regularity conditions ensuring among other things that the convergence
in (\ref{Convergence: Ghat}) and (\ref{Convergence: Q}) is suitably uniform in
$\mathbf{s,}$ (\ref{Cube root asymptotics}) then follows from an application
of a continuous mapping-type theorem for the $\operatorname*{argmax}$
functional to the representation in (\ref{Representation: Thetahat}).

Similarly to (\ref{Representation: Thetahat}), the bootstrap analogue of
$\mathbf{\hat{\theta}}_{n}$ admits a representation of the form%
\[
\sqrt[3]{n}(\mathbf{\hat{\theta}}_{n}^{\ast}-\mathbf{\hat{\theta}}%
_{n})=\operatorname*{argmax}\limits_{\mathbf{s}\in\mathbb{R}^{d}}\{\hat{G}%
_{n}^{\ast}(\mathbf{s})+\hat{Q}_{n}(\mathbf{s})\},
\]
where, for $\mathbf{s}$ such that $\mathbf{\hat{\theta}}_{n}+\mathbf{s}%
n^{-1/3}\in\mathbf{\Theta}$,%
\[
\hat{G}_{n}^{\ast}(\mathbf{s})=n^{2/3}[\hat{M}_{n}^{\ast}(\mathbf{\hat{\theta
}}_{n}+\mathbf{s}n^{-1/3})-\hat{M}_{n}^{\ast}(\mathbf{\hat{\theta}}_{n}%
)-\hat{M}_{n}(\mathbf{\hat{\theta}}_{n}+\mathbf{s}n^{-1/3})+\hat{M}%
_{n}(\mathbf{\hat{\theta}}_{n})]
\]
and%
\[
\hat{Q}_{n}(\mathbf{s})=n^{2/3}[\hat{M}_{n}(\mathbf{\hat{\theta}}%
_{n}+\mathbf{s}n^{-1/3})-\hat{M}_{n}(\mathbf{\hat{\theta}}_{n})].
\]
Under mild conditions, $\hat{G}_{n}^{\ast}$ satisfies the following bootstrap
counterpart of (\ref{Convergence: Ghat}):%
\begin{equation}
\hat{G}_{n}^{\ast}(\mathbf{s})\rightsquigarrow_{\mathbb{P}}\mathcal{G}%
_{0}(\mathbf{s}),\label{Convergence: Ghatstar}%
\end{equation}
where $\rightsquigarrow_{\mathbb{P}}$ denotes conditional weak convergence in
probability
%TCIMACRO{\TeXButton{\citep*[defined as in][Section
%2.9]{vanderVaart-Wellner_1996_Book}}{\citep
%*[defined as in][Section 2.9]{vanderVaart-Wellner_1996_Book}}}%
%BeginExpansion
\citep*[defined as in][Section 2.9]{vanderVaart-Wellner_1996_Book}%
%EndExpansion
. On the other hand, although $\hat{Q}_{n}$ is non-random under the bootstrap
distribution and satisfies $\operatorname*{argmax}_{\mathbf{s}\in
\mathbb{R}^{d}}\hat{Q}_{n}(\mathbf{s})=\mathbf{0,}$ it turns out that $\hat
{Q}_{n}(\mathbf{s})\nrightarrow_{\mathbb{P}}\mathcal{Q}_{0}(\mathbf{s})$ in
general. In other words, the natural bootstrap counterpart of
(\ref{Convergence: Q}) typically fails and, as a partial consequence, so does
the natural bootstrap counterpart of (\ref{Cube root asymptotics}); i.e.,
$\sqrt[3]{n}(\mathbf{\hat{\theta}}_{n}^{\ast}-\mathbf{\hat{\theta}}%
_{n})\not \rightsquigarrow _{\mathbb{P}}\operatorname*{argmax}_{\mathbf{s}%
\in\mathbb{R}^{d}}\{\mathcal{G}_{0}(\mathbf{s})+\mathcal{Q}_{0}(\mathbf{s})\}$.

To the extent that the inconsistency of the bootstrap can be attributed to the
fact that the shape of $\hat{Q}_{n}$ fails to replicate that of $Q_{n},$ it
seems plausible that a consistent bootstrap-based distributional approximation
can be obtained by basing the approximation on%
\[
\mathbf{\tilde{\theta}}_{n}^{\ast}=\operatorname*{argmax}%
\limits_{\mathbf{\theta\in\Theta}}\tilde{M}_{n}^{\ast}(\mathbf{\theta}%
),\qquad\tilde{M}_{n}^{\ast}(\mathbf{\theta})=\frac{1}{n}\sum_{i=1}^{n}%
\tilde{m}_{n}(\mathbf{z}_{i,n}^{\ast},\mathbf{\theta}),
\]
where $\tilde{m}_{n}$ is a suitably \textquotedblleft
reshaped\textquotedblright\ version of $m_{0}$ satisfying two properties.
First, $\tilde{G}_{n}^{\ast}$ should be asymptotically equivalent to $\hat
{G}_{n}^{\ast},$ where $\tilde{G}_{n}^{\ast}$ is the counterpart of $\hat
{G}_{n}^{\ast}$ associated with $\tilde{m}_{n}:$%
\[
\tilde{G}_{n}^{\ast}(\mathbf{s})=n^{2/3}[\tilde{M}_{n}^{\ast}(\mathbf{\hat
{\theta}}_{n}+\mathbf{s}n^{-1/3})-\tilde{M}_{n}^{\ast}(\mathbf{\hat{\theta}%
}_{n})-\tilde{M}_{n}(\mathbf{\hat{\theta}}_{n}+\mathbf{s}n^{-1/3})+\tilde
{M}_{n}(\mathbf{\hat{\theta}}_{n})],\qquad\tilde{M}_{n}(\mathbf{\theta}%
)=\frac{1}{n}\sum_{i=1}^{n}\tilde{m}_{n}(\mathbf{z}_{i},\mathbf{\theta}).
\]
Second, and most importantly, $\tilde{Q}_{n}$ should be asymptotically
quadratic, where $\tilde{Q}_{n}$ is the counterpart of $\hat{Q}_{n}$
associated with $\tilde{m}_{n}$:%
\[
\tilde{Q}_{n}(\mathbf{s})=n^{2/3}[\tilde{M}_{n}(\mathbf{\hat{\theta}}%
_{n}+\mathbf{s}n^{-1/3})-\tilde{M}_{n}(\mathbf{\hat{\theta}}_{n})].
\]

Accordingly, let%
\[
\tilde{m}_{n}(\mathbf{z},\mathbf{\theta})=m_{0}(\mathbf{z},\mathbf{\theta
})-\hat{M}_{n}(\mathbf{\theta})-\frac{1}{2}(\mathbf{\theta}-\mathbf{\hat
{\theta}}_{n})^{\prime}\mathbf{\tilde{H}}_{n}(\mathbf{\theta}-\mathbf{\hat
{\theta}}_{n}),
\]
where $\mathbf{\tilde{H}}_{n}$ is an estimator of $\mathbf{H}_{0}.$ Then%
\[
\sqrt[3]{n}(\mathbf{\tilde{\theta}}_{n}^{\ast}-\mathbf{\hat{\theta}}%
_{n})=\operatorname*{argmax}\limits_{\mathbf{s}\in\mathbb{R}^{d}}\{\tilde
{G}_{n}^{\ast}(\mathbf{s})+\tilde{Q}_{n}(\mathbf{s})\},
\]
where, by construction, $\tilde{G}_{n}^{\ast}(\mathbf{s})=\hat{G}_{n}^{\ast
}(\mathbf{s})$ and $\tilde{Q}_{n}(\mathbf{s})=-\mathbf{s}^{\prime
}\mathbf{\tilde{H}}_{n}\mathbf{s}/2\mathbf{.}$ Because $\tilde{G}_{n}^{\ast
}=\hat{G}_{n}^{\ast},$ $\tilde{G}_{n}^{\ast}(\mathbf{s})\rightsquigarrow
_{\mathbb{P}}\mathcal{G}_{0}(\mathbf{s})$ whenever
(\ref{Convergence: Ghatstar}) holds. In addition, $\tilde{Q}_{n}%
(\mathbf{s})\rightarrow_{\mathbb{P}}\mathcal{Q}_{0}(\mathbf{s})$ provided
$\mathbf{\tilde{H}}_{n}\rightarrow_{\mathbb{P}}\mathbf{H}_{0}.$ As a
consequence, it seems plausible that if $\mathbf{\tilde{H}}_{n}\rightarrow
_{\mathbb{P}}\mathbf{H}_{0},$ then $\sqrt[3]{n}(\mathbf{\tilde{\theta}}%
_{n}^{\ast}-\mathbf{\hat{\theta}}_{n})\rightsquigarrow_{\mathbb{P}%
}\operatorname*{argmax}\limits_{\mathbf{s}\in\mathbb{R}^{d}}\{\mathcal{G}%
_{0}(\mathbf{s})+\mathcal{Q}_{0}(\mathbf{s})\}.$

For the purposes of situating this paper in the literature, the following
alternative heuristic explanation of our approach may be useful. Restating the
result in (\ref{Cube root asymptotics}) as%
\[
\sqrt[3]{n}(\mathbf{\hat{\theta}}_{n}-\mathbf{\theta}_{0})\rightsquigarrow
\mathcal{S}_{0}(\mathcal{G}_{0}),\qquad\mathcal{S}_{0}(\mathcal{G}%
)=\operatorname*{argmax}\limits_{\mathbf{s}\in\mathbb{R}^{d}}\{\mathcal{G}%
(\mathbf{s})+\mathcal{Q}_{0}(\mathbf{s})\},
\]
our procedure approximates the distribution of $\mathcal{S}_{0}(\mathcal{G}%
_{0})$ by that of $\mathcal{\tilde{S}}_{n}(\hat{G}_{n}^{\ast}),$ where the
distribution of the bootstrap process $\hat{G}_{n}^{\ast}\ $approximates that
of $\mathcal{G}_{0}$ and where $\mathcal{\tilde{S}}_{n}(\mathcal{G}%
)=\operatorname*{argmax}\limits_{\mathbf{s}\in\mathbb{R}^{d}}\{\mathcal{G}%
(\mathbf{s})+\tilde{Q}_{n}(\mathbf{s})\}$ is an estimator of $\mathcal{S}%
_{0}(\mathcal{G}).$ In other words, our procedure replaces the functional
$\mathcal{S}_{0}$ with a consistent estimator (namely, $\mathcal{\tilde{S}%
}_{n}$) and its random argument $\mathcal{G}_{0}$ with a bootstrap
approximation (namely, $\hat{G}_{n}^{\ast}$). The same type of generic
construction has appeared in the econometrics literature before, notably in
\cite{Andrews-Soares_2010_ECMA} and \cite{Fang-Santos_2019_REStud}.

Our bootstrap-based distributional approximation can be shown to be consistent
also in the more standard case where $m_{n}(\mathbf{z},\mathbf{\theta})$ is
sufficiently smooth in $\mathbf{\theta}$ to ensure that an approximate
maximizer of $\hat{M}_{n}$ is asymptotically normal and that the nonparametric
bootstrap is consistent. In fact, $\mathbf{\tilde{\theta}}_{n}^{\ast}$ is
(first-order) asymptotically equivalent to $\mathbf{\hat{\theta}}_{n}^{\ast}$
in that standard case, so our procedure can be interpreted as a modification
of the nonparametric bootstrap that is designed to be \textquotedblleft
robust\textquotedblright\ to the types of non-smoothness that give rise to
cube root asymptotics.

\section{Main Result\label{[Section] Main Result}}

When making the heuristics of Section \ref{[Section] Heuristics} precise, we
consider the more general situation where the estimator $\mathbf{\hat{\theta}%
}_{n}$ is an approximate maximizer (with respect to $\mathbf{\theta}%
\in\mathbf{\Theta}\subseteq\mathbb{R}^{d}$) of%
\[
\hat{M}_{n}(\mathbf{\theta})=\frac{1}{n}\sum_{i=1}^{n}m_{n}(\mathbf{z}%
_{i},\mathbf{\theta}),
\]
where $m_{n}$ is a known function, and where $\mathbf{z}_{1},\ldots
,\mathbf{z}_{n}$ is a random sample of a random vector $\mathbf{z}.$ This
formulation of $\hat{M}_{n},$ which reduces to that of Section
\ref{[Section] Heuristics} when $m_{n}$ does not depend on $n,$ is adopted in
order to cover certain estimation problems where, rather than admitting a
characterization of the form (\ref{Estimand}), the estimand $\mathbf{\theta
}_{0}$ admits the characterization%
\[
\mathbf{\theta}_{0}=\operatorname*{argmax}\limits_{\mathbf{\theta}%
\in\mathbf{\Theta}}M_{0}(\mathbf{\theta}),\qquad M_{0}(\mathbf{\theta}%
)=\lim_{n\rightarrow\infty}M_{n}(\mathbf{\theta}),\qquad M_{n}(\mathbf{\theta
})=\mathbb{E}[m_{n}(\mathbf{z},\mathbf{\theta})].
\]

In other words, in the setting considered in this section, $\mathbf{\hat
{\theta}}_{n}$ approximately maximizes a function $\hat{M}_{n}$ whose
population counterpart $M_{n}$ can be interpreted as a regularization
%TCIMACRO{\TeXButton{\citep*[in the sense of][]{Bickel-Li_2006_Test}}%
%{\citep*[in the sense of][]{Bickel-Li_2006_Test}} }%
%BeginExpansion
\citep*[in the sense of][]{Bickel-Li_2006_Test}
%EndExpansion
of a function $M_{0}$ whose maximizer $\mathbf{\theta}_{0}$ is the object of
interest. This generalization is attractive because it allows us to formulate
results that cover local $M$-estimators such as the conditional maximum score
estimator of \cite{Honore-Kyriazidou_2000_ECMA}. Studying this setting,
\cite{Seo-Otsu_2018_AoS} gave conditions under which $\mathbf{\hat{\theta}%
}_{n}$ converges at a rate equal to the cube root of the \textquotedblleft
effective\textquotedblright\ sample size and has a limiting distribution of
\cite{Chernoff_1964_AISM} type. Analogous conclusions will be drawn below,
albeit under slightly different conditions.

For any $n$ and any $\delta>0,$ define%
\[
\bar{m}_{n}(\mathbf{z})=\sup_{m\in\mathcal{M}_{n}}|m(\mathbf{z})|,\qquad
\mathcal{M}_{n}=\{m_{n}(\cdot,\mathbf{\theta}):\mathbf{\theta}\in
\mathbf{\Theta}\},
\]
and%
\[
\bar{d}_{n}^{\delta}(\mathbf{z})=\sup_{d\in\mathcal{D}_{n}^{\delta}%
}|d(\mathbf{z})|,\qquad\mathcal{D}_{n}^{\delta}=\{m_{n}(\cdot,\mathbf{\theta
})-m_{n}(\cdot,\mathbf{\theta}_{0}):\mathbf{\theta}\in\mathbf{\Theta}%
_{0}^{\delta}\},\qquad\mathbf{\Theta}_{0}^{\delta}=\{\mathbf{\theta}%
\in\mathbf{\Theta}:||\mathbf{\theta-\theta}_{0}||\leq\delta\}.
\]

\begin{description}
\item[Condition CRA (Cube Root Asymptotics)] For some $q_{n}>0$ with
$r_{n}=\sqrt[3]{nq_{n}}\rightarrow\infty,$ the following are
satisfied:\newline(i) $\{\mathcal{M}_{n}:n\geq1\}$ is uniformly manageable for
the envelopes $\bar{m}_{n}$ and $q_{n}\mathbb{E}[\bar{m}_{n}(\mathbf{z}%
)^{2}]=O(1).$\newline Also, $\sup_{\mathbf{\theta\in\Theta}}\left\vert
M_{n}(\mathbf{\theta})-M_{0}(\mathbf{\theta})\right\vert =o(1)$ and, for every
$\delta>0,$ $\sup_{\mathbf{\theta}\in\mathbf{\Theta}\backslash\mathbf{\Theta
}_{0}^{\delta}}M_{0}(\mathbf{\theta})<M_{0}(\mathbf{\theta}_{0}).$\newline(ii)
$\mathbf{\theta}_{0}$ is an interior point of $\mathbf{\Theta}$ and, for some
$\delta>0,$ $M_{0}$ and $M_{n}$ are twice continuously differentiable on
$\mathbf{\Theta}_{0}^{\delta}$ and $\sup_{\mathbf{\theta}\in\mathbf{\Theta
}_{0}^{\delta}}\left\Vert \partial^{2}[M_{n}(\mathbf{\theta})-M_{0}%
(\mathbf{\theta})]/\partial\mathbf{\theta}\partial\mathbf{\theta}^{\prime
}\right\Vert =o(1).$\newline Also, $r_{n}||\partial M_{n}(\mathbf{\theta}%
_{0})/\partial\mathbf{\theta}||=o(1)$ and $\mathbf{H}_{0}=-\partial^{2}%
M_{0}(\mathbf{\theta}_{0})/\partial\mathbf{\theta}\partial\mathbf{\theta
}^{\prime}$ is positive definite.\newline(iii) For some $\delta>0,$
$\{\mathcal{D}_{n}^{\delta^{\prime}}:n\geq1,0<\delta^{\prime}\leq\delta\}$ is
uniformly manageable for the envelopes $\bar{d}_{n}^{\delta^{\prime}}$ and
$q_{n}\sup_{0<\delta^{\prime}\leq\delta}\mathbb{E}[\bar{d}_{n}^{\delta
^{\prime}}(\mathbf{z})^{2}/\delta^{\prime}]=O(1).$\newline(iv) For every
$\delta_{n}>0$ with $\delta_{n}=O(r_{n}^{-1}),$ $q_{n}^{3}r_{n}^{-1}%
\mathbb{E}[\bar{d}_{n}^{\delta_{n}}(\mathbf{z})^{4}]=o(1) $ and, for all
$\mathbf{s},\mathbf{t}\in\mathbb{R}^{d}$ and for some $\mathcal{C}_{0}$ with
$\mathcal{C}_{0}(\mathbf{s},\mathbf{s})+\mathcal{C}_{0}(\mathbf{t}%
,\mathbf{t})-2\mathcal{C}_{0}(\mathbf{s},\mathbf{t})>0$ for $\mathbf{s\neq
t,}$%
\[
\sup_{\mathbf{\theta}\in\mathbf{\Theta}_{0}^{\delta_{n}}}\left\vert
\frac{q_{n}}{\delta_{n}}\mathbb{E}[\{m_{n}(\mathbf{z},\mathbf{\theta}%
+\delta_{n}\mathbf{s})-m_{n}(\mathbf{z},\mathbf{\theta})\}\{m_{n}%
(\mathbf{z},\mathbf{\theta}+\delta_{n}\mathbf{t})-m_{n}(\mathbf{z}%
,\mathbf{\theta})\}]-\mathcal{C}_{0}(\mathbf{s},\mathbf{t})\right\vert =o(1).
\]
(v) For every $\delta_{n}>0$ with $\delta_{n}=O(r_{n}^{-1}),$%
\[
\lim_{C\rightarrow\infty}\underset{n\rightarrow\infty}{\lim\sup}\sup
_{0<\delta\leq\delta_{n}}q_{n}\mathbb{E}[%
%TCIMACRO{\TeXButton{I}{\I}}%
%BeginExpansion
\I
%EndExpansion
\{q_{n}\bar{d}_{n}^{\delta}(\mathbf{z})>C\}\bar{d}_{n}^{\delta}(\mathbf{z}%
)^{2}/\delta]=0
\]
and $\sup_{\mathbf{\theta,\theta}^{\prime}\in\mathbf{\Theta}_{0}^{\delta_{n}}%
}\mathbb{E}[|m_{n}(\mathbf{z},\mathbf{\theta})-m_{n}(\mathbf{z},\mathbf{\theta
}^{\prime})|]/||\mathbf{\theta-\theta}^{\prime}||=O(1).$
\end{description}

To interpret Condition CRA, consider first the benchmark case where
$m_{n}=m_{0}$ and $q_{n}=1.$ In this case, the condition is similar to
assumptions (ii)-(vii) of the main theorem of \cite{Kim-Pollard_1990_AoS}, to
which the reader is referred for a definition of the term (uniformly)
manageable. The differences between their assumptions and Condition CRA are
technical in nature, since we need to slightly strengthen their assumptions in
order to be able to analyze the bootstrap. For instance, the displayed part of
Condition CRA(iv) is a locally uniform (with respect to $\mathbf{\theta}$ near
$\mathbf{\theta}_{0}$) version of its counterpart in
\cite{Kim-Pollard_1990_AoS}. More generally, Condition CRA can be interpreted
as an $n$-varying version of a suitably (for the purpose of analyzing the
bootstrap) strengthened version of the assumptions of
\cite{Kim-Pollard_1990_AoS}. The differences between Condition CRA and the
$i.i.d. $ version of the conditions in \cite{Seo-Otsu_2018_AoS} are also
technical in nature, but for completeness we highlight two here. First, they
control the complexity of various function classes using the concept of
bracketing entropy, while we follow \cite{Kim-Pollard_1990_AoS} and obtain
maximal inequalities using bounds on uniform entropy numbers implied by the
concept of (uniform) manageability. Second, whereas \cite{Seo-Otsu_2018_AoS}
control the bias of $\mathbf{\hat{\theta}}_{n}$ through a locally uniform
bound on $M_{n}-M_{0},$ Condition CRA controls the bias through the first and
second derivatives of $M_{n}-M_{0}.$

Under Condition CRA, the effective sample size is $nq_{n}=r_{n}^{3}$ and if
$\mathbf{\hat{\theta}}_{n}$ is an approximate maximizer of $\hat{M}_{n},$ then
$r_{n}(\mathbf{\hat{\theta}}_{n}-\mathbf{\theta}_{0})$ has a limiting
distribution of \cite{Chernoff_1964_AISM} type. The heuristics of the previous
section are rate-adaptive (i.e., $\sqrt[3]{n}$ can be replaced by a generic
$r_{n}$), so once again it stands to reason that if $\mathbf{\tilde{H}}_{n}$
is a consistent estimator of $\mathbf{H}_{0},$ then the distribution of
$r_{n}(\mathbf{\hat{\theta}}_{n}-\mathbf{\theta}_{0})$ can be consistently
estimated by that of $r_{n}(\mathbf{\tilde{\theta}}_{n}^{\ast}-\mathbf{\hat
{\theta}}_{n}),$ where $\mathbf{\tilde{\theta}}_{n}^{\ast}$ is an approximate
maximizer of%
\[
\tilde{M}_{n}^{\ast}(\mathbf{\theta})=\frac{1}{n}\sum_{i=1}^{n}\tilde{m}%
_{n}(\mathbf{z}_{i,n}^{\ast},\mathbf{\theta}),\qquad\tilde{m}_{n}%
(\mathbf{z},\mathbf{\theta})=m_{n}(\mathbf{z},\mathbf{\theta})-\hat{M}%
_{n}(\mathbf{\theta})-\frac{1}{2}(\mathbf{\theta}-\mathbf{\hat{\theta}}%
_{n})^{\prime}\mathbf{\tilde{H}}_{n}(\mathbf{\theta}-\mathbf{\hat{\theta}}%
_{n}),
\]
with $\mathbf{z}_{1,n}^{\ast},\ldots,\mathbf{z}_{n,n}^{\ast}$ being a random
sample from the empirical distribution of $\mathbf{z}_{1},\ldots
,\mathbf{z}_{n}$. A precise statement is given in the following theorem.

\begin{theorem}
\label{[Theorem] Main Result}Suppose Condition CRA holds. If $\mathbf{\tilde
{H}}_{n}\rightarrow_{\mathbb{P}}\mathbf{H}_{0}$ and if%
\[
\hat{M}_{n}(\mathbf{\hat{\theta}}_{n})\geq\sup_{\mathbf{\theta}\in
\mathbf{\Theta}}\hat{M}_{n}(\mathbf{\theta})-o_{\mathbb{P}}(r_{n}%
^{-2})\text{\qquad and\qquad}\tilde{M}_{n}^{\ast}(\mathbf{\tilde{\theta}}%
_{n}^{\ast})\geq\sup_{\mathbf{\theta}\in\mathbf{\Theta}}\tilde{M}_{n}^{\ast
}(\mathbf{\theta})-o_{\mathbb{P}}(r_{n}^{-2}),
\]
then%
\begin{equation}
r_{n}(\mathbf{\hat{\theta}}_{n}-\mathbf{\theta}_{0})\rightsquigarrow
\operatorname*{argmax}\limits_{\mathbf{s}\in\mathbb{R}^{d}}\{\mathcal{G}%
_{0}(\mathbf{s})+\mathcal{Q}_{0}(\mathbf{s}%
)\},\label{Nonparametric cube root asymptotics}%
\end{equation}
and%
\begin{equation}
r_{n}(\mathbf{\tilde{\theta}}_{n}^{\ast}-\mathbf{\hat{\theta}}_{n}%
)\rightsquigarrow_{\mathbb{P}}\operatorname*{argmax}_{\mathbf{s}\in
\mathbb{R}^{d}}\{\mathcal{G}_{0}(\mathbf{s})+\mathcal{Q}_{0}(\mathbf{s}%
)\},\label{Nonparametric cube root asymptotics (bootstrap)}%
\end{equation}
where $\mathcal{G}_{0}$ is a zero-mean Gaussian process with covariance kernel
$\mathcal{C}_{0}$ and $\mathcal{Q}_{0}(\mathbf{s})=-\mathbf{s}^{\prime
}\mathbf{H}_{0}\mathbf{s}/2.$
\end{theorem}

The algorithm for our proposed bootstrap-based distributional approximation is
as follows:\newline\emph{Step 1}. Using the sample $\mathbf{z}_{1}%
,\dots,\mathbf{z}_{n},$ compute $\mathbf{\hat{\theta}}_{n}$ by approximately
maximizing $\hat{M}_{n}(\mathbf{\theta}).$\newline\emph{Step 2}. Using
$\mathbf{\hat{\theta}}_{n}$ and $\mathbf{z}_{1},\dots,\mathbf{z}_{n},$ compute
$\mathbf{\tilde{H}}_{n}.$ (A generic estimator $\mathbf{\tilde{H}}_{n}$ is
given in Section \ref{[Subsection] Estimation of H0}.)\newline\emph{Step 3}.
Using $\mathbf{\hat{\theta}}_{n},$ $\mathbf{\tilde{H}}_{n},$ and the bootstrap
sample $\mathbf{z}_{1,n}^{\ast},\dots,\mathbf{z}_{n,n}^{\ast},$ compute
$\mathbf{\tilde{\theta}}_{n}^{\ast}$ by approximately maximizing $\tilde
{M}_{n}^{\ast}(\mathbf{\theta}).$ ($\mathbf{\hat{\theta}}_{n}$ and
$\mathbf{\tilde{H}}_{n}$ are not recomputed at this step.)\newline\emph{Step
4}. Repeat \emph{Step 3} to generate draws from the distribution of
$r_{n}(\mathbf{\tilde{\theta}}_{n}^{\ast}-\mathbf{\hat{\theta}}_{n}).$

\subsection{Estimation of $\mathbf{H}_{0}$%
\label{[Subsection] Estimation of H0}}

A generic numerical derivative estimator of $\mathbf{H}_{0}$ is the matrix
$\mathbf{\tilde{H}}_{n}^{\mathtt{ND}}$ with element $(k,l)$ given by%
\[
\tilde{H}_{n,kl}^{\mathtt{ND}}=-\frac{1}{4\epsilon_{n}^{2}}[\hat{M}%
_{n}(\mathbf{\hat{\theta}}_{n}+\mathbf{e}_{k}\epsilon_{n}+\mathbf{e}%
_{l}\epsilon_{n})-\hat{M}_{n}(\mathbf{\hat{\theta}}_{n}+\mathbf{e}_{k}%
\epsilon_{n}-\mathbf{e}_{l}\epsilon_{n})-\hat{M}_{n}(\mathbf{\hat{\theta}}%
_{n}-\mathbf{e}_{k}\epsilon_{n}+\mathbf{e}_{l}\epsilon_{n})+\hat{M}%
_{n}(\mathbf{\hat{\theta}}_{n}-\mathbf{e}_{k}\epsilon_{n}-\mathbf{e}%
_{l}\epsilon_{n})],
\]
where $\mathbf{e}_{k}$ is the $k$th unit vector in $\mathbb{R}^{d}$ and where
$\epsilon_{n}$ is a positive tuning parameter. Conditions under which this
estimator is consistent are given in the following lemma.

\begin{lemma}
\label{[Lemma] Consistency of VtildeND}Suppose Condition CRA holds and that
$r_{n}(\mathbf{\hat{\theta}}_{n}-\mathbf{\theta}_{0})=O_{\mathbb{P}}(1).$ If
$\epsilon_{n}\rightarrow0$ and if $r_{n}\epsilon_{n}\rightarrow\infty,$ then
$\mathbf{\tilde{H}}_{n}^{\mathtt{ND}}\rightarrow_{\mathbb{P}}\mathbf{H}_{0}.$
\end{lemma}

Plausibility of the high-level condition $r_{n}(\mathbf{\hat{\theta}}%
_{n}-\mathbf{\theta}_{0})=O_{\mathbb{P}}(1)$ follows from
(\ref{Nonparametric cube root asymptotics}). To facilitate practical
implementation, it is useful to go beyond consistency and develop a Nagar-type
mean squared error (MSE) expansion that can be used to select $\epsilon_{n}.$
To state one such result for $\tilde{H}_{n,kl}^{\mathtt{ND}},$ define%
\[
\ddot{M}_{n,kl}(\mathbf{\theta})=\frac{\partial^{2}}{\partial\theta
_{k}\partial\theta_{l}}M_{n}(\mathbf{\theta}),\qquad\ddot{M}_{0,kl}%
(\mathbf{\theta})=\frac{\partial^{2}}{\partial\theta_{k}\partial\theta_{l}%
}M_{0}(\mathbf{\theta}),
\]%
\[
\mathsf{B}_{kl}=-\frac{1}{6}\left[  \frac{\partial^{2}}{\partial\theta_{k}%
^{2}}\ddot{M}_{0,kl}(\mathbf{\theta}_{0})+\frac{\partial^{2}}{\partial
\theta_{l}^{2}}\ddot{M}_{0,kl}(\mathbf{\theta}_{0})\right]  ,
\]
and%
\[
\mathsf{V}_{kl}=\frac{1}{8}[\mathcal{C}_{0}(\mathbf{e}_{k}+\mathbf{e}%
_{l},\mathbf{e}_{k}+\mathbf{e}_{l})+\mathcal{C}_{0}(\mathbf{e}_{k}%
-\mathbf{e}_{l},\mathbf{e}_{k}-\mathbf{e}_{l})-2\mathcal{C}_{0}(\mathbf{e}%
_{k}+\mathbf{e}_{l},\mathbf{e}_{k}-\mathbf{e}_{l})-2\mathcal{C}_{0}%
(\mathbf{e}_{k}+\mathbf{e}_{l},-\mathbf{e}_{k}+\mathbf{e}_{l})].
\]

\begin{lemma}
\label{[Lemma] MSE of VtildeND}Suppose the conditions of Lemma
\ref{[Lemma] Consistency of VtildeND} hold and that, for some $\delta>0,$
$\ddot{M}_{0,kl}$ and $\ddot{M}_{n,kl}$ are twice continuously differentiable
on $\mathbf{\Theta}_{0}^{\delta}$ with $\sup_{\mathbf{\theta}\in
\mathbf{\Theta}_{0}^{\delta}}||\partial^{2}[\ddot{M}_{n,kl}(\mathbf{\theta
})-\ddot{M}_{0,kl}(\mathbf{\theta})]/\partial\mathbf{\theta}\partial
\mathbf{\theta}^{\prime}||=o(1).$ If $\mathcal{C}_{0}(\mathbf{s}%
,-\mathbf{s})=0$ and $\mathcal{C}_{0}(\mathbf{s},\mathbf{t})=\mathcal{C}%
_{0}(-\mathbf{s},-\mathbf{t})$ for all $\mathbf{s},\mathbf{t}\in\mathbb{R}%
^{d},$ then $\tilde{H}_{n,kl}^{\mathtt{ND}}$ admits an approximation
$\check{H}_{n,kl}^{\mathtt{ND}}$ satisfying
\[
\tilde{H}_{n,kl}^{\mathtt{ND}}=\check{H}_{n,kl}^{\mathtt{ND}}+o_{\mathbb{P}%
}\left(  \epsilon_{n}^{2}+\frac{1}{\sqrt{r_{n}^{3}\epsilon_{n}^{3}}}\right)
+O_{\mathbb{P}}\left(  \frac{1}{r_{n}}\right)  ,
\]
where the $O_{\mathbb{P}}(1/r_{n})$ term does not depend on $\epsilon_{n}$ and
where%
\[
\mathbb{E}[(\check{H}_{n,kl}^{\mathtt{ND}}-H_{n,kl})^{2}]=\epsilon_{n}%
^{4}\mathsf{B}_{kl}^{2}+\frac{1}{r_{n}^{3}\epsilon_{n}^{3}}\mathsf{V}%
_{kl}+o\left(  \epsilon_{n}^{4}+\frac{1}{r_{n}^{3}\epsilon_{n}^{3}}\right)
,\qquad H_{n,kl}=-\ddot{M}_{n,kl}(\mathbf{\theta}_{0}).
\]

\end{lemma}

The conditions $\mathcal{C}_{0}(\mathbf{s},-\mathbf{s})=0$ and $\mathcal{C}%
_{0}(\mathbf{s},\mathbf{t})=\mathcal{C}_{0}(-\mathbf{s},-\mathbf{t})$ are
satisfied in all of the examples we have analyzed. Using the lemma, the
approximate MSE (AMSE), $\epsilon_{n}^{4}\mathsf{B}_{kl}^{2}+r_{n}%
^{-3}\epsilon_{n}^{-3}\mathsf{V}_{kl}$, can be minimized by choosing
$\epsilon_{n}$ proportional to $r_{n}^{-3/7},$ the optimal factor of
proportionality being a function of $\mathsf{B}_{kl}$ and $\mathsf{V}_{kl}.$
To be specific, the optimal $\epsilon_{n}$ is given by $\epsilon
_{n,kl}^{\mathtt{AMSE}}=(3\mathsf{V}_{kl}/4\mathsf{B}_{kl}^{2})^{1/7}%
r_{n}^{-3/7},$ a feasible version of which can be constructed by replacing
$\mathsf{B}_{kl}$ and $\mathsf{V}_{kl}$ with preliminary estimators thereof.

\section{Examples\label{[Section] Examples}}

\subsection{Maximum Score\label{[Subsection] Maximum Score}}

To describe a version of the maximum score estimator of
%TCIMACRO{\TeXButton{\cite{Manski_1975_JoE,Manski_1985_JoE}}{\cite
%{Manski_1975_JoE,Manski_1985_JoE}}}%
%BeginExpansion
\cite{Manski_1975_JoE,Manski_1985_JoE}%
%EndExpansion
, suppose $\mathbf{z}_{1},\ldots,\mathbf{z}_{n}$ is a random sample of
$\mathbf{z}=(y,\mathbf{x}^{\prime})^{\prime}$ generated by the binary response
model%
\[
y=%
%TCIMACRO{\TeXButton{I}{\I}}%
%BeginExpansion
\I
%EndExpansion
(\mathbf{x}^{\prime}\mathbf{\beta}_{0}+u\geq0),\qquad\mathrm{Median}%
(u|\mathbf{x})=0,
\]
where $\mathbf{\beta}_{0}\in\mathbb{R}^{d+1}$ is an unknown parameter of
interest, $\mathbf{x}\in\mathbb{R}^{d+1}$ is a vector of covariates, and $u$
is an unobserved error term. Following \cite{Abrevaya-Huang_2005_ECMA}, we
employ the parameterization $\mathbf{\beta}_{0}=(1,\mathbf{\theta}_{0}%
^{\prime})^{\prime},$ where $\mathbf{\theta}_{0}\in\mathbb{R}^{d}$ is unknown.
In other words, we assume that the first element of $\mathbf{\beta}_{0}$ is
positive and then normalize the (unidentified) scale of $\mathbf{\beta}_{0}$
by setting its first element equal to unity. Partitioning $\mathbf{x}$
conformably with $\mathbf{\beta}_{0}$ as $\mathbf{x}=(x_{1},\mathbf{x}%
_{2}^{\prime})^{\prime},$ a maximum score estimator of $\mathbf{\theta}_{0}$
is any $\mathbf{\hat{\theta}}_{n}^{\mathtt{MS}}$ approximately maximizing
$\hat{M}_{n}$ for $m_{n}(\mathbf{z},\mathbf{\theta})=m^{\mathtt{MS}%
}(\mathbf{z},\mathbf{\theta})=(2y-1)%
%TCIMACRO{\TeXButton{I}{\I}}%
%BeginExpansion
\I
%EndExpansion
(x_{1}+\mathbf{x}_{2}^{\prime}\mathbf{\theta}\geq0)$.

Regarded as a member of the class of $M$-estimators exhibiting cube root
asymptotics, the maximum score estimator is representative in a couple of
respects. First, under easy-to-interpret primitive conditions the estimator is
covered by the results of Section \ref{[Section] Main Result}. Second, in
addition to the generic estimator $\mathbf{\tilde{H}}_{n}^{\mathtt{ND}}$
discussed above, the maximum score estimator admits example-specific
consistent estimators of the $\mathbf{H}_{0}$ associated with it.

Under standard regularity conditions (stated in Section A.2 of the
supplemental appendix), Condition CRA is satisfied with $q_{n}=1,$%
\[
\mathbf{H}_{0}=\mathbf{H}^{\mathtt{MS}}=2\mathbb{E}[f_{u|x_{1},\mathbf{x}_{2}%
}(0|-\mathbf{x}_{2}^{\prime}\mathbf{\theta}_{0},\mathbf{x}_{2})f_{x_{1}%
|\mathbf{x}_{2}}(-\mathbf{x}_{2}^{\prime}\mathbf{\theta}_{0}|\mathbf{x}%
_{2})\mathbf{x}_{2}\mathbf{x}_{2}^{\prime}],
\]
and%
\[
\mathcal{C}_{0}(\mathbf{s},\mathbf{t})=\mathcal{C}^{\mathtt{MS}}%
(\mathbf{s},\mathbf{t})=\mathbb{E}[f_{x_{1}|\mathbf{x}_{2}}(-\mathbf{x}%
_{2}^{\prime}\mathbf{\theta}_{0}|\mathbf{x}_{2})\min\{|\mathbf{x}_{2}^{\prime
}\mathbf{s}|,|\mathbf{x}_{2}^{\prime}\mathbf{t}|\}%
%TCIMACRO{\TeXButton{I}{\I}}%
%BeginExpansion
\I
%EndExpansion
(\operatorname*{sgn}(\mathbf{x}_{2}^{\prime}\mathbf{s})=\operatorname*{sgn}%
(\mathbf{x}_{2}^{\prime}\mathbf{t}))],
\]
where $f_{a|\mathbf{b}}$ denotes the conditional Lebesgue density of $a$ given
$\mathbf{b}.$ As a consequence, Theorem \ref{[Theorem] Main Result} is
applicable to $\mathbf{\hat{\theta}}_{n}^{\mathtt{MS}}$ and the consistency
requirement $\mathbf{\tilde{H}}_{n}\rightarrow_{\mathbb{P}}\mathbf{H}%
^{\mathtt{MS}}$ is satisfied by the numerical derivative estimator discussed
in Section \ref{[Subsection] Estimation of H0} if $\epsilon_{n}\rightarrow0 $
and $n\epsilon_{n}^{3}\rightarrow\infty.$ Under the additional regularity
conditions of Lemma \ref{[Lemma] MSE of VtildeND}, MSE-optimal tuning
parameter choices are feasible. In addition, alternative consistent estimators
of $\mathbf{H}^{\mathtt{MS}}$ can be constructed exploiting the specific
structure of this example. One option is to employ a \textquotedblleft
plug-in\textquotedblright\ estimator of $\mathbf{H}^{\mathtt{MS}}$ based on
nonparametric estimators of $f_{u|x_{1},\mathbf{x}_{2}}$ and $f_{x_{1}%
|\mathbf{x}_{2}}.$ An alternative, example-specific estimator is%
\[
\mathbf{\tilde{H}}_{n}^{\mathtt{MS}}=-\frac{1}{n}\sum_{i=1}^{n}(2y_{i}%
-1)\dot{K}_{n}(x_{1i}+\mathbf{x}_{2i}^{\prime}\mathbf{\hat{\theta}}%
_{n}^{\mathtt{MS}})\mathbf{x}_{2i}\mathbf{x}_{2i}^{\prime},
\]
where, for a differentiable kernel function $K$ and a positive bandwidth
$h_{n},$ $\dot{K}_{n}(u)=dK_{n}(u)/du$ and $K_{n}(u)=K(u/h_{n})/h_{n}.$ In
words, $\mathbf{\tilde{H}}_{n}^{\mathtt{MS}}$ is constructed by
\textquotedblleft smoothing out\textquotedblright\ the indicator function
entering $m^{\mathtt{MS}}(\mathbf{z},\mathbf{\theta})$ and then twice
differentiating the corresponding\ objective function
%TCIMACRO{\TeXButton{\citep*[previously used by][]{Horowitz_1992_ECMA}}%
%{\citep*[previously used by][]{Horowitz_1992_ECMA}}}%
%BeginExpansion
\citep*[previously used by][]{Horowitz_1992_ECMA}%
%EndExpansion
.

\subsection{Panel Maximum Score\label{[Subsection] Panel Maximum Score}}

Consider the panel data binary response model%
\[
Y_{t}=%
%TCIMACRO{\TeXButton{I}{\I}}%
%BeginExpansion
\I
%EndExpansion
(\mathbf{X}_{t}^{\prime}\mathbf{\beta}_{0}+\alpha+u_{t}\geq0),\qquad t=1,2,
\]
where $\mathbf{\beta}_{0}\in\mathbb{R}^{d+1}$ is an unknown parameter of
interest, $\alpha$ is an unobserved (time-invariant) individual-specific
effect, and $u_{t}$ is an unobserved error term. Analyzing this model,
\cite{Manski_1987_ECMA} gave conditions under which $\mathbf{\beta}_{0}$ is
identified up to scale and demonstrated consistency of a conditional maximum
score estimator.

Suppose $\mathbf{\beta}_{0}$ is identified up to scale and that its first
element is positive, in which case we can normalize that element to unity and
employ the parameterization $\mathbf{\beta}_{0}=(1,\mathbf{\theta}_{0}%
^{\prime})^{\prime},$ where $\mathbf{\theta}_{0}\in\mathbb{R}^{d}$ is unknown.
To describe a version of the estimator of \cite{Manski_1987_ECMA}, partition
$\mathbf{X}_{t}$ conformably with $\mathbf{\beta}_{0}$ as $\mathbf{X}%
_{t}=(X_{1t},\mathbf{X}_{2t}^{\prime})^{\prime}$ and define $\mathbf{z}%
=(y,x_{1},\mathbf{x}_{2}^{\prime})^{\prime},$ where $y=Y_{2}-Y_{1},$
$x_{1}=X_{12}-X_{11},$ and $\mathbf{x}_{2}=(\mathbf{X}_{22}-\mathbf{X}_{21}).$
Assuming $\mathbf{z}_{1},\ldots,\mathbf{z}_{n}$ is a random sample of
$\mathbf{z},$ a panel maximum score estimator of $\mathbf{\theta}_{0}\ $is any
$\mathbf{\hat{\theta}}_{n}^{\mathtt{PMS}}$ approximately maximizing $\hat
{M}_{n}$ for $m_{n}(\mathbf{z},\mathbf{\theta})=m^{\mathtt{PMS}}%
(\mathbf{z},\mathbf{\theta})=y%
%TCIMACRO{\TeXButton{I}{\I}}%
%BeginExpansion
\I
%EndExpansion
(x_{1}+\mathbf{x}_{2}^{\prime}\mathbf{\theta}\geq0).$

As one would expect, the properties of $\mathbf{\hat{\theta}}_{n}%
^{\mathtt{PMS}}$ are qualitatively similar to those of $\mathbf{\hat{\theta}%
}_{n}^{\mathtt{MS}}.$ To be specific, under regularity conditions (stated in
Section A.3 of the supplemental appendix), the panel maximum score estimator
is covered by the results of Section \ref{[Section] Main Result} and an
example-specific alternative to the generic numerical derivative estimator is
available, namely%
\[
\mathbf{\tilde{H}}_{n}^{\mathtt{PMS}}=-n^{-1}\sum_{i=1}^{n}y_{i}\dot{K}%
_{n}(x_{1i}+\mathbf{x}_{2i}^{\prime}\mathbf{\hat{\theta}}_{n}^{\mathtt{PMS}%
})\mathbf{x}_{2i}\mathbf{x}_{2i}^{\prime},
\]
where $\dot{K}_{n}$ is as in the maximum score example.

\subsection{Conditional Maximum
Score\label{[Subsection] Conditional Maximum Score}}

Consider the dynamic panel data binary response model%
\[
Y_{t}=%
%TCIMACRO{\TeXButton{I}{\I}}%
%BeginExpansion
\I
%EndExpansion
(\mathbf{X}_{t}^{\prime}\mathbf{\beta}_{0}+Y_{t-1}\gamma_{0}+\alpha+u_{t}%
\geq0),\qquad t=1,2,3,
\]
where $\mathbf{\beta}_{0}\in\mathbb{R}^{d}$ and $\gamma_{0}\in\mathbb{R}$ are
unknown parameters of interest, $\alpha$ is an unobserved (time-invariant)
individual-specific effect, and $u_{t}$ is an unobserved error term.
\cite{Honore-Kyriazidou_2000_ECMA} analyzed this model and gave conditions
under which $\mathbf{\beta}_{0}$ and $\gamma_{0}$ are identified up to a
common scale factor. Assuming these conditions hold and that the first element
of $\mathbf{\beta}_{0}$ is positive, we can normalize that element to unity
and employ the parameterization $(\mathbf{\beta}_{0}^{\prime},\gamma
_{0})^{\prime}=(1,\mathbf{\theta}_{0}^{\prime})^{\prime},$ where
$\mathbf{\theta}_{0}\in\mathbb{R}^{d}$ is unknown.

To describe a version of the conditional maximum score estimator of
\cite{Honore-Kyriazidou_2000_ECMA}, partition $\mathbf{X}_{t}$ after the first
element as $\mathbf{X}_{t}=(X_{1t},\mathbf{X}_{2t}^{\prime})^{\prime}$ and
define $\mathbf{z}=(y,x_{1},\mathbf{x}_{2}^{\prime},\mathbf{w}^{\prime
})^{\prime},$ where $y=Y_{2}-Y_{1},$ $x_{1}=X_{12}-X_{11},$ $\mathbf{x}%
_{2}=((\mathbf{X}_{22}-\mathbf{X}_{21})^{\prime},Y_{3}-Y_{0})^{\prime},$ and
$\mathbf{w}=\mathbf{X}_{2}-\mathbf{X}_{3}.$ Assuming $\mathbf{z}_{1}%
,\ldots,\mathbf{z}_{n}$ is a random sample of $\mathbf{z},$ a conditional
maximum score estimator of $\mathbf{\theta}_{0}$ is any $\mathbf{\hat{\theta}%
}_{n}^{\mathtt{CMS}}$ approximately maximizing $\hat{M}_{n}$ for
$m_{n}(\mathbf{z},\mathbf{\theta})=m_{n}^{\mathtt{CMS}}(\mathbf{z}%
,\mathbf{\theta})=y%
%TCIMACRO{\TeXButton{I}{\I}}%
%BeginExpansion
\I
%EndExpansion
(x_{1}+\mathbf{x}_{2}^{\prime}\mathbf{\theta}\geq0)\kappa_{n}(\mathbf{w}),$
where $\kappa_{n}(\mathbf{w})=\kappa(\mathbf{w}/b_{n})/b_{n}^{d}$ for a kernel
function $\kappa$ and a bandwidth $b_{n}.$

Through its dependence on $b_{n},$ the function $m_{n}^{\mathtt{CMS}}$ depends
on $n$ in a non-negligible way. In particular, the effective sample size is
$nb_{n}^{d}$ (rather than $n$) in the current setting, so to the extent that
they exist one would expect primitive sufficient conditions for Condition CRA
to include $q_{n}=b_{n}^{d}$ in this example. Apart from this predictable
change, the properties of the conditional maximum score estimator
$\mathbf{\hat{\theta}}_{n}^{\mathtt{CMS}}$ turn out to be qualitatively
similar to those of $\mathbf{\hat{\theta}}_{n}^{\mathtt{MS}}$. To be specific,
under regularity conditions (stated in Section A.4 of the supplemental
appendix), the conditional maximum score estimator is covered by the results
of Section \ref{[Section] Main Result} and an example-specific alternative to
the generic numerical derivative estimator is available, namely%
\[
\mathbf{\tilde{H}}_{n}^{\mathtt{CMS}}=-n^{-1}\sum_{i=1}^{n}y_{i}\dot{K}%
_{n}(x_{1i}+\mathbf{x}_{2i}^{\prime}\mathbf{\hat{\theta}}_{n}^{\mathtt{CMS}%
})\mathbf{x}_{2i}\mathbf{x}_{2i}^{\prime}\kappa_{n}(\mathbf{w}_{i}),
\]
where $\dot{K}_{n}\ $is as in the maximum score example.

\subsection{Empirical Risk
Minimization\label{[Subsection] Empirical Risk Minimization}}

\cite{Mohammadi-vandeGeer_2005_JMLR} considered two-category classification
problems in machine learning. Specifically, given a binary outcome
$y\in\{-1,1\}$ and a vector of features $\mathbf{x}\in\mathcal{X},$ the goal
is to estimate the $\mathbf{\theta}_{0}$ that minimizes the misclassification
error (or risk) $\mathbb{P}[h_{\mathbf{\theta}}(\mathbf{x})\neq y]$ with
respect to $\mathbf{\theta}\in\mathbf{\Theta}\subseteq\mathbb{R}^{d},$ where
$\{h_{\mathbf{\theta}}:\mathbf{\theta}\in\mathbf{\Theta}\}$ is a collection of
classifiers. For simplicity, we consider the case where the feature is
univariate with support $\mathcal{X}=[0,1]$ and the classifiers are of the
form%
\[
h_{\mathbf{\theta}}(x)=\sum_{\ell=1}^{d+1}(-1)^{\ell}%
%TCIMACRO{\TeXButton{I}{\I}}%
%BeginExpansion
\I
%EndExpansion
(\theta_{\ell-1}\leq x<\theta_{\ell}),\qquad\mathbf{\theta}=(\theta_{1}%
,\theta_{2},\cdots,\theta_{d})^{\prime},
\]
where $\mathbf{\Theta}=\{(\theta_{1},\theta_{2},\cdots,\theta_{d})^{\prime}%
\in\lbrack0,1]^{d}:0=\theta_{0}\leq\theta_{1}\leq\cdots\leq\theta_{d}%
\leq\theta_{d+1}=1\}.$

Assuming $\mathbf{z}_{1},\ldots,\mathbf{z}_{n}$ is a random sample of
$\mathbf{z},$ an empirical risk minimizer is any $\mathbf{\hat{\theta}}%
_{n}^{\mathtt{ERM}}$ approximately maximizing $\hat{M}_{n}$ for $m_{n}%
(\mathbf{z},\mathbf{\theta})=m^{\mathtt{ERM}}(\mathbf{z},\mathbf{\theta})=-%
%TCIMACRO{\TeXButton{I}{\I}}%
%BeginExpansion
\I
%EndExpansion
(h_{\mathbf{\theta}}(x)\neq y).$ Under regularity conditions similar to those
of \cite[Section 2.1]{Mohammadi-vandeGeer_2005_JMLR}, the empirical risk
minimizer is covered by Theorem \ref{[Theorem] Main Result} and the
consistency requirement on $\mathbf{\tilde{H}}_{n}$ can be met in various
ways; for details, see Section A.5 of the supplemental appendix.

\section{Simulations\label{[Section] Simulations}}

We illustrate the numerical performance of our proposed bootstrap-based
inference methods for the maximum score estimator. Given the setup in Section
\ref{[Subsection] Maximum Score}, we generate data from that model with $d=1,$
$\mathbf{\theta}_{0}=1,$ $\mathbf{x}=(x_{1},x_{2})^{\prime}\thicksim
\mathcal{N}((0,1)^{\prime},\mathbf{I}_{2})$ with $\mathbf{I}_{2}$ the
$(2\times2)$ identity matrix, and $u$ generated by three distinct
distributions. Specifically, DGP 1 sets $u\thicksim\mathsf{Logistic}%
(0,1)/\sqrt{2\pi^{2}/3},$ DGP 2 sets $u\thicksim\mathcal{T}_{3}/\sqrt{3}$,
where $\mathcal{T}_{3}$ denotes a Student's $t$-distribution with 3 degrees of
freedom, and DGP 3 sets $u\thicksim(1+2(x_{1}+x_{2})^{2}+(x_{1}+x_{2}%
)^{4})\mathsf{Logistic}(0,1)/\sqrt{\pi^{2}/48}.$

The Monte Carlo experiment employs a sample size $n=1,000$ with $B=2,000$
bootstrap replications and $S=2,000$ simulations. For each of the three DGPs,
we implement the standard non-parametric bootstrap, the $m$-out-of-$n$
bootstrap using $m\in\{\left\lceil n^{1/2}\right\rceil ,\left\lceil
n^{2/3}\right\rceil ,\left\lceil n^{4/5}\right\rceil \},$ and our proposed
method using the two estimators $\mathbf{\tilde{H}}_{n}^{\mathtt{MS}}$ and
$\mathbf{\tilde{H}}_{n}^{\mathtt{ND}}$ of $\mathbf{H}_{0}.$ We report
empirical coverage for nominal $95\%$ confidence intervals and their average
interval length. For the case of our proposed procedures, we investigate their
performance using (i)\ infeasible (simulation-based)\ MSE-optimal choices of
tuning parameters (bandwidth/derivative step), denoted by $h_{\mathtt{MSE}}$
and $\epsilon_{\mathtt{MSE}}$, and (ii)\ infeasible and feasible AMSE-optimal
choices of the tuning parameters, denoted by $h_{\mathtt{AMSE}}$, $\hat
{h}_{\mathtt{AMSE}}$, $\epsilon_{\mathtt{AMSE}}$ and $\hat{\epsilon
}_{\mathtt{AMSE}}$; for details, see Section A.2 of the supplemental appendix.

\singlespacing
\renewcommand{\arraystretch}{1.1}
\newcommand{\pathsims}{PUT-LAST-WHAT-YOU-WANT}
\renewcommand{\pathsims}%
{C:/Documents/Research/Cattaneo-Jansson-Nagasawa--BootstrapCubeRoot}
\renewcommand{\pathsims}{Z:/research/Cattaneo-Jansson-Nagasawa_2020_ECMA}
\clearpage\begin{table}\renewcommand{\arraystretch}{1.2}
\caption{Simulations, Maximum Score Estimator, 95\% Confidence Intervals.}%
\label{table:maxscore}
\resizebox{\textwidth}{!}{%latex.default(round(out, 3), file = "table_maxscore.txt", append = FALSE,     table.env = FALSE, center = "none", title = "", n.cgroup = c(3,         3, 3), cgroup = c("DGP 1", "DGP 2", "DGP 3"), colheads = rep(c("$h,\\epsilon$",         "Coverage", "Length"), 3), n.rgroup = c(1, 3, 3, 3),     rgroup = c("Standard", "m-out-of-n", "Plug-in: $\\tilde{\\mathbf{V}}^{\\mathtt{MS}}_n$",         "Num Deriv: $\\tilde{\\mathbf{V}}^{\\mathtt{ND}}_n$"),     rowname = c("", "$m = \\lceil n^{1/2} \\rceil$", "$m = \\lceil n^{2/3} \\rceil$",         "$m = \\lceil n^{4/5} \\rceil$", "$h_{ \\mathtt{MSE} } $",         "$h_{\\mathtt{AMSE} }$", "$\\hat{h}_{\\mathtt{AMSE} }$",         "$\\epsilon_{ \\mathtt{MSE} } $", "$\\epsilon_{\\mathtt{AMSE} }$",         "$\\hat{\\epsilon}_{\\mathtt{AMSE} }$"))%
\begin{tabular}{lrrrcrrrcrrr}
\hline\hline
\multicolumn{1}{l}{\bfseries }&\multicolumn{3}{c}{\bfseries DGP 1}&\multicolumn{1}{c}{\bfseries }&\multicolumn{3}{c}{\bfseries DGP 2}&\multicolumn{1}{c}{\bfseries }&\multicolumn{3}{c}{\bfseries DGP 3}\tabularnewline
\cline{2-4} \cline{6-8} \cline{10-12}
\multicolumn{1}{l}{}&\multicolumn{1}{c}{$h,\epsilon$}&\multicolumn{1}{c}{Coverage}&\multicolumn{1}{c}{Length}&\multicolumn{1}{c}{}&\multicolumn{1}{c}{$h,\epsilon$}&\multicolumn{1}{c}{Coverage}&\multicolumn{1}{c}{Length}&\multicolumn{1}{c}{}&\multicolumn{1}{c}{$h,\epsilon$}&\multicolumn{1}{c}{Coverage}&\multicolumn{1}{c}{Length}\tabularnewline
\hline
{\bfseries Standard}&&&&&&&&&&&\tabularnewline
~~&$$&$0.625$&$0.472$&&$$&$0.647$&$0.475$&&$$&$0.654$&$0.243$\tabularnewline
\hline
{\bfseries m-out-of-n}&&&&&&&&&&&\tabularnewline
~~$m = \lceil n^{1/2} \rceil$&$$&$0.997$&$1.698$&&$$&$0.998$&$1.753$&&$$&$1.000$&$1.890$\tabularnewline
~~$m = \lceil n^{2/3} \rceil$&$$&$0.978$&$1.185$&&$$&$0.983$&$1.221$&&$$&$0.989$&$0.724$\tabularnewline
~~$m = \lceil n^{4/5} \rceil$&$$&$0.899$&$0.820$&&$$&$0.897$&$0.837$&&$$&$0.930$&$0.447$\tabularnewline
\hline
{\bfseries Plug-in: $\tilde{\mathbf{V}}^{\mathtt{MS}}_n$}&&&&&&&&&&&\tabularnewline
~~$h_{ \mathtt{MSE} } $&$0.620$&$0.954$&$0.511$&&$0.580$&$0.957$&$0.523$&&$0.150$&$0.962$&$0.277$\tabularnewline
~~$h_{\mathtt{AMSE} }$&$1.108$&$0.972$&$0.590$&&$0.480$&$0.951$&$0.518$&&$0.123$&$0.942$&$0.263$\tabularnewline
~~$\hat{h}_{\mathtt{AMSE} }$&$0.443$&$0.940$&$0.508$&&$0.409$&$0.946$&$0.518$&&$0.155$&$0.957$&$0.278$\tabularnewline
\hline
{\bfseries Num Deriv: $\tilde{\mathbf{V}}^{\mathtt{ND}}_n$}&&&&&&&&&&&\tabularnewline
~~$\epsilon_{ \mathtt{MSE} } $&$1.400$&$0.936$&$0.483$&&$1.360$&$0.938$&$0.485$&&$0.290$&$0.939$&$0.249$\tabularnewline
~~$\epsilon_{\mathtt{AMSE} }$&$0.537$&$0.880$&$0.414$&&$0.573$&$0.894$&$0.426$&&$0.224$&$0.902$&$0.227$\tabularnewline
~~$\hat{\epsilon}_{\mathtt{AMSE} }$&$0.518$&$0.876$&$0.413$&&$0.512$&$0.882$&$0.420$&&$0.369$&$0.947$&$0.270$\tabularnewline
\hline
\end{tabular}
}
\footnotesize\hspace{.5in} Notes:\newline(i) Panel \textbf{Standard}
refers to standard nonparametric bootstrap, Panel \textbf{m-out-of-n}
refers to $m$-out-of-$n$ nonparametric bootstrap with subsample size $m$, Panel \textbf
{Plug-in:} $\tilde{\mathbf{V}}^\mathtt{MS}%
_n$ refers to our proposed bootstrap-based implemented using the example-specific plug-in drift estimator, and Panel \textbf
{Num Deriv:} $\tilde{\mathbf{V}}^\mathtt{ND}%
_n$ refers to our proposed bootstrap-based implemented using the generic numerical derivative drift estimator.\newline
(ii) Column ``$h$, $\epsilon
$'' reports tuning parameter value used or average across simulations when estimated, and Columns ``Coverage'' and ``Length'' report empirical coverage and average length of bootstrap-based $95\%$ percentile confidence intervals, respectively.\newline
(iii) $h_\mathtt{MSE}$ and $\epsilon_\mathtt{MSE}%
$ correspond to the simulation MSE-optimal choices, $h_\mathtt{AMSE}%
$ and $\epsilon_\mathtt{AMSE}%
$ correspond to the AMSE-optimal choices, and $\hat{h}_\mathtt{AMSE}%
$ and $\hat{\epsilon}_\mathtt{AMSE}%
$ correspond to the ROT feasible implementation of $\hat{h}_\mathtt
{AMSE}$ and $\hat{\epsilon}_\mathtt{AMSE}%
$ described in the supplemental appendix.
\end{table}%
%EndExpansion

\doublespacing

Table
%TCIMACRO{\TeXButton{\ref{table:maxscore}}{\ref{table:maxscore}} }%
%BeginExpansion
\ref{table:maxscore}
%EndExpansion
presents the main results, which are consistent across all three simulation
designs. First, as expected, the standard nonparametric bootstrap (labeled
\textquotedblleft Standard\textquotedblright)\ does not perform well, leading
to confidence intervals with an average $64\%$ empirical coverage rate.
Second, the $m$-out-of-$n$ bootstrap (labeled \textquotedblleft
m-out-of-n\textquotedblright) performs somewhat better for small subsamples,
but leads to very large average interval length of the resulting confidence
intervals. Our proposed methods, on the other hand, exhibit good finite sample
performance in this Monte Carlo experiment. Results employing the
example-specific plug-in estimator $\mathbf{\tilde{H}}_{n}^{\mathtt{MS}}$ are
presented under the label \textquotedblleft Plug-in\textquotedblright\ while
results employing the generic numerical derivative estimator $\mathbf{\tilde
{H}}_{n}^{\mathtt{ND}}$ are reported under the label \textquotedblleft Num
Deriv\textquotedblright. Empirical coverage appears stable across different
values of the tuning parameters for each method, with better performance in
the case of $\mathbf{\tilde{H}}_{n}^{\mathtt{MS}}$. We conjecture that
$n=1,000$ is too small for the numerical derivative estimator $\mathbf{\tilde
{H}}_{n}^{\mathtt{ND}}$ to lead to as good inference performance as
$\mathbf{\tilde{H}}_{n}^{\mathtt{MS}}$ (e.g., note that the MSE-optimal choice
$\epsilon_{\mathtt{MSE}}$ is greater than $1$). Nevertheless, empirical
coverage of confidence intervals constructed using our proposed
bootstrap-based method is close to $95\%$ in all cases except when
$\mathbf{\tilde{H}}_{n}^{\mathtt{ND}}$ is used with either the infeasible
asymptotic choice $\epsilon_{\mathtt{AMSE}}$ or its estimated counterpart
$\hat{\epsilon}_{\mathtt{AMSE}},$ and with an average interval length of at
most half that of any of the $m$-out-of-$n$ competing confidence intervals. In
particular, confidence intervals based on $\mathbf{\tilde{H}}_{n}%
^{\mathtt{MS}}$ implemented with the feasible bandwidth $\hat{h}%
_{\mathtt{AMSE}}$ perform quite well across the three DGPs considered.

\section{Conclusion\label{[Section] Conclusion}}

We developed a valid resampling procedure for cube root asymptotics based on
the nonparametric bootstrap. Whereas the standard nonparametric bootstrap is
known to be invalid in the setting we study, we show that bootstrap validity
can be restored by applying a carefully tailored reshapement of the objective
function defining the estimator. Such reshapement is easy to implement both in
general and in specific cases, as illustrated by the distinct examples we considered.

\cite{Seo-Otsu_2018_AoS} gave conditions under which results of the form
(\ref{Nonparametric cube root asymptotics}) can be obtained also when the data
exhibits weak dependence; see also \cite{Bagchi-Banerjee-Stoev_2016_JASA}, and
references therein. It seems plausible that a version of our procedure,
implemented with a resampling procedure suitable for dependent data, can be
shown to be consistent under similar conditions, but it is beyond the scope of
this paper to substantiate that conjecture.

\section{Proof of Theorem \ref{[Theorem] Main Result}%
\label{[Section] Proof of Theorem 1}}

The proof proceeds by first showing (\ref{Nonparametric cube root asymptotics}%
) and then using that result to establish
(\ref{Nonparametric cube root asymptotics (bootstrap)}). In both cases, we
employ arguments similar to those used in the proof of the main theorem of
\cite{Kim-Pollard_1990_AoS}. The remainder of this section outlines the main
steps in the proof; for technical details, see Lemmas A.1-A.10 in Section A.1
of the supplemental appendix.

\textit{Proof of (\ref{Nonparametric cube root asymptotics}).} The estimator
$\mathbf{\hat{\theta}}_{n}$ is assumed to satisfy%
\[
\left.  \{\hat{G}_{n}(\mathbf{s})+Q_{n}(\mathbf{s})\}\right\vert
_{\mathbf{s=}r_{n}(\mathbf{\hat{\theta}}_{n}-\mathbf{\theta}_{0})}\geq
\sup_{\mathbf{s}\in\mathbb{R}^{d}}\{\hat{G}_{n}(\mathbf{s})+Q_{n}%
(\mathbf{s})\}+o_{\mathbb{P}}(1),
\]
where%
\[
\hat{G}_{n}(\mathbf{s})=r_{n}^{2}[\hat{M}_{n}(\mathbf{\theta}_{0}%
+\mathbf{s}r_{n}^{-1})-\hat{M}_{n}(\mathbf{\theta}_{0})-M_{n}(\mathbf{\theta
}_{0}+\mathbf{s}r_{n}^{-1})+M_{n}(\mathbf{\theta}_{0})]%
%TCIMACRO{\TeXButton{I}{\I}}%
%BeginExpansion
\I
%EndExpansion
(\mathbf{\theta}_{0}+\mathbf{s}r_{n}^{-1}\in\mathbf{\Theta})
\]
and%
\[
Q_{n}(\mathbf{s})=r_{n}^{2}[M_{n}(\mathbf{\theta}_{0}+\mathbf{s}r_{n}%
^{-1})-M_{n}(\mathbf{\theta}_{0})]%
%TCIMACRO{\TeXButton{I}{\I}}%
%BeginExpansion
\I
%EndExpansion
(\mathbf{\theta}_{0}+\mathbf{s}r_{n}^{-1}\in\mathbf{\Theta}).
\]
By the argmax continuous mapping theorem
%TCIMACRO{\TeXButton{\citep*[e.g.,][Theorem
%3.2.2]{vanderVaart-Wellner_1996_Book}}{\citep
%*[e.g.,][Theorem 3.2.2]{vanderVaart-Wellner_1996_Book}}}%
%BeginExpansion
\citep*[e.g.,][Theorem 3.2.2]{vanderVaart-Wellner_1996_Book}%
%EndExpansion
, it therefore suffices to show that $r_{n}(\mathbf{\hat{\theta}}%
_{n}-\mathbf{\theta}_{0})=O_{\mathbb{P}}(1)$ and that $\hat{G}_{n}%
+Q_{n}\rightsquigarrow\mathcal{G}_{0}+\mathcal{Q}_{0}$ in the topology of
uniform convergence on compacta. (The other conditions required by the argmax
continuous mapping theorem are easily verified.)

To obtain the rate of convergence of $\mathbf{\hat{\theta}}_{n},$ we begin by
using a standard argument to show that $\mathbf{\hat{\theta}}_{n}%
-\mathbf{\theta}_{0}=o_{\mathbb{P}}(1)$ under Condition CRA(i) and then
strengthen that conclusion to $r_{n}(\mathbf{\hat{\theta}}_{n}-\mathbf{\theta
}_{0})=O_{\mathbb{P}}(1)$ by using Conditions CRA(ii)-(iii) and proceeding
along the lines of \cite[Theorem 3.2.5]{vanderVaart-Wellner_1996_Book}. In
both cases, we employ the maximal inequality in \cite[Theorem 4.2]%
{Pollard_1989_SS}; for details, see Lemmas A.1 and A.3 of the supplemental appendix.

Next, because $Q_{n}$ is non-random, $\hat{G}_{n}+Q_{n}\rightsquigarrow
\mathcal{G}_{0}+\mathcal{Q}_{0}$ in the topology of uniform convergence on
compacta if $Q_{n}$ converges compactly to $\mathcal{Q}_{0}$ and if $\hat
{G}_{n}\rightsquigarrow\mathcal{G}_{0}$ in the topology of uniform convergence
on compacta. Compact convergence of $Q_{n}$ follows from Condition CRA (ii);
for details, see Lemma A.2 of the supplemental appendix. Also, to show that
$\hat{G}_{n}\rightsquigarrow\mathcal{G}_{0}$ in the topology of uniform
convergence on compacta, it suffices to show that $\hat{G}_{n}$ converges to
$\mathcal{G}_{0}$ in the sense of weak convergence of finite-dimensional
projections and that $\{\hat{G}_{n}(\mathbf{s}):||\mathbf{s}||\leq K\}$ is
stochastically equicontinuous for every $K>0.$

Under Conditions CRA(ii)-(iv), weak convergence of finite-dimensional
projections can be shown using the Cram\'{e}r-Wold device and the fact that
$\mathbb{E}[\hat{G}_{n}(\mathbf{s})\hat{G}_{n}(\mathbf{t}))]$ converges to
$\mathcal{C}_{0}(\mathbf{s},\mathbf{t})$ for every $\mathbf{s},\mathbf{t}%
\in\mathbb{R}^{d};$ for details, see Lemma A.4 of the supplemental appendix.
Finally, under Conditions CRA(iii) and CRA(v) and employing the maximal
inequality in \cite[Theorem 4.2]{Pollard_1989_SS}, stochastic equicontinuity
can be shown by proceeding as in the proof of \cite[Lemma 4.6]%
{Kim-Pollard_1990_AoS}; for details, see Lemma A.5 of the supplemental appendix.

\textit{Proof of (\ref{Nonparametric cube root asymptotics (bootstrap)}).} The
proof of (\ref{Nonparametric cube root asymptotics (bootstrap)}) is a natural
bootstrap analog of the proof of (\ref{Nonparametric cube root asymptotics}).
The estimator $\mathbf{\tilde{\theta}}_{n}^{\ast}$ is assumed to satisfy%
\[
\left.  \{\tilde{G}_{n}^{\ast}(\mathbf{s})+\tilde{Q}_{n}(\mathbf{s}%
)\}\right\vert _{\mathbf{s=}r_{n}(\mathbf{\tilde{\theta}}_{n}^{\ast
}-\mathbf{\hat{\theta}}_{n})}\geq\sup_{\mathbf{s}\in\mathbb{R}^{d}}\{\tilde
{G}_{n}^{\ast}(\mathbf{s})+\tilde{Q}_{n}(\mathbf{s})\}+o_{\mathbb{P}}(1),
\]
where%
\[
\tilde{G}_{n}^{\ast}(\mathbf{s})=r_{n}^{2}[\tilde{M}_{n}^{\ast}(\mathbf{\hat
{\theta}}_{n}+\mathbf{s}r_{n}^{-1})-\tilde{M}_{n}^{\ast}(\mathbf{\hat{\theta}%
}_{n})-\tilde{M}_{n}(\mathbf{\hat{\theta}}_{n}+\mathbf{s}r_{n}^{-1})+\tilde
{M}_{n}^{\ast}(\mathbf{\hat{\theta}}_{n})]%
%TCIMACRO{\TeXButton{I}{\I}}%
%BeginExpansion
\I
%EndExpansion
(\mathbf{\hat{\theta}}_{n}+\mathbf{s}r_{n}^{-1}\in\mathbf{\Theta})
\]
and%
\[
\tilde{Q}_{n}(\mathbf{s})=r_{n}^{2}[\tilde{M}_{n}(\mathbf{\hat{\theta}}%
_{n}+\mathbf{s}r_{n}^{-1})-\tilde{M}_{n}(\mathbf{\hat{\theta}}_{n})]%
%TCIMACRO{\TeXButton{I}{\I}}%
%BeginExpansion
\I
%EndExpansion
(\mathbf{\hat{\theta}}_{n}+\mathbf{s}r_{n}^{-1}\in\mathbf{\Theta})=-\frac
{1}{2}\mathbf{s}^{\prime}\mathbf{\tilde{H}}_{n}\mathbf{s}%
%TCIMACRO{\TeXButton{I}{\I}}%
%BeginExpansion
\I
%EndExpansion
(\mathbf{\hat{\theta}}_{n}+\mathbf{s}r_{n}^{-1}\in\mathbf{\Theta}).
\]
By the argmax continuous mapping theorem, it therefore suffices to show that
$r_{n}(\mathbf{\tilde{\theta}}_{n}^{\ast}-\mathbf{\hat{\theta}}_{n}%
)=O_{\mathbb{P}}(1)$ and that $\tilde{G}_{n}^{\ast}+\tilde{Q}_{n}%
\rightsquigarrow_{\mathbb{P}}\mathcal{G}_{0}+\mathcal{Q}_{0}$ in the topology
of uniform convergence on compacta.

Using $\mathbf{\tilde{H}}_{n}\rightarrow_{\mathbb{P}}\mathbf{H}_{0},$ to
obtain the rate of convergence of $\mathbf{\tilde{\theta}}_{n}^{\ast}$ we
first show that $\mathbf{\tilde{\theta}}_{n}^{\ast}-\mathbf{\hat{\theta}}%
_{n}=o_{\mathbb{P}}(1)$ under Condition CRA(i) and then strengthen that
conclusion to $r_{n}(\mathbf{\tilde{\theta}}_{n}^{\ast}-\mathbf{\hat{\theta}%
}_{n})=O_{\mathbb{P}}(1)$ by using $r_{n}(\mathbf{\hat{\theta}}_{n}%
-\mathbf{\theta}_{0})=O_{\mathbb{P}}(1)$ and Condition CRA(iii). As in the
derivation of the convergence rate of $\mathbf{\hat{\theta}}_{n},$ both steps
employ the maximal inequality in \cite[Theorem 4.2]{Pollard_1989_SS}; for
details, see Lemmas A.6 and A.8 of the supplemental appendix.

Next, because $\mathcal{Q}_{0}$ is non-random, $\tilde{G}_{n}^{\ast}+\tilde
{Q}_{n}\rightsquigarrow_{\mathbb{P}}\mathcal{G}_{0}+\mathcal{Q}_{0}$ in the
topology of uniform convergence on compacta if $\tilde{Q}_{n}\rightarrow
_{\mathbb{P}}\mathcal{Q}_{0}$ in the topology of uniform convergence on
compacta and if $\hat{G}_{n}\rightsquigarrow\mathcal{G}_{0}$ in the topology
of uniform convergence on compacta. By construction, $\tilde{Q}_{n}$ is such
that if $\mathbf{\tilde{H}}_{n}\rightarrow_{\mathbb{P}}\mathbf{H}_{0}$ and if
$\mathbf{\hat{\theta}}_{n}\rightarrow_{\mathbb{P}}\mathbf{\theta}_{0}%
\in\operatorname*{int}(\mathbf{\Theta),}$ then $\tilde{Q}_{n}\rightarrow
_{\mathbb{P}}\mathcal{Q}_{0}$ in the topology of uniform convergence on
compacta; for details, see Lemma A.7 of the supplemental appendix.

Also, to show that $\tilde{G}_{n}^{\ast}\rightsquigarrow_{\mathbb{P}%
}\mathcal{G}_{0}$ in the topology of uniform convergence on compacta, it
suffices to show that $\tilde{G}_{n}$ converges to $\mathcal{G}_{0}$ in the
sense of conditional weak convergence in probability of finite-dimensional
projections and that $\{\tilde{G}_{n}^{\ast}(\mathbf{s}):||\mathbf{s}||\leq
K\}$ is stochastically equicontinuous for every $K>0.$ Conditional weak
convergence in probability of finite-dimensional projections can be shown
using the Cram\'{e}r-Wold device and the fact that the maximal inequality in
\cite[Theorem 4.2]{Pollard_1989_SS} can be used to show that $\mathbb{E}%
_{n}^{\ast}[\tilde{G}_{n}^{\ast}(\mathbf{s})\tilde{G}_{n}^{\ast}%
(\mathbf{t}))]$ converges in probability to $\mathcal{C}_{0}(\mathbf{s}%
,\mathbf{t})$ for every $\mathbf{s},\mathbf{t}\in\mathbb{R}^{d}$, where
$\mathbb{E}_{n}^{\ast}$ denotes an expectation computed under the bootstrap
distribution conditional on the data; for details, see Lemma A.9 of the
supplemental appendix. Finally, employing the maximal inequality in
\cite[Theorem 4.2]{Pollard_1989_SS}, stochastic equicontinuity can be shown by
proceeding as in the proof of \cite[Lemma 4.6]{Kim-Pollard_1990_AoS}; for
details, see Lemma A.10 of the supplemental appendix.

\singlespacing

\bibliographystyle{econometrica}
\bibliography{BootstrapCubeRoot}

\end{document}